\documentclass[smallextended,envcountsect]{svjour3}
\usepackage{graphicx,color,cite}
\usepackage{mathrsfs}
\usepackage{amsmath, amscd, amsfonts, amssymb, graphicx, color}
\usepackage[bookmarksnumbered, colorlinks, plainpages]{hyperref}
\usepackage{mathptmx}   
\usepackage[colorinlistoftodos]{todonotes}
\usepackage[left=4cm, right=4cm, top=3.5cm, bottom=3.5cm]{geometry}

\usepackage{latexsym}

\if{
\newtheorem{theorem}{Theorem}[section]
\newtheorem{definition}{Definition}[section]

\newtheorem{remark}{Remark}[section]
\newtheorem{example}{Example}[section]
\newtheorem{claim}{\textbf{Claim}}[section]

\newtheorem{corollary}{Corollary}[section]

\newtheorem{claim}{\textbf{Claim}}[section]

}\fi

\newcommand{\ball}{\mathbb{B}}

\renewcommand\theenumi{\roman{enumi}}

\newcounter{mycount}

\smartqed

\usepackage{etex}

\usepackage{mathtools}

\usepackage{enumitem}
\renewcommand\theenumi{(\roman{enumi})}
\renewcommand\theenumii{(\alph{enumii})}

\usepackage{csquotes}

\usepackage{todonotes}

\makeatletter
\let\orgdescriptionlabel\descriptionlabel
\renewcommand*{\descriptionlabel}[1]{
 \let\orglabel\label
 \let\label\@gobble
 \phantomsection
 \edef\@currentlabel{#1}
 \let\label\orglabel
 \orgdescriptionlabel{#1}
}
\def\th@plain{
 \thm@notefont{}
 \itshape
}
\def\th@definition{
 \thm@notefont{}
 \normalfont
}

\g@addto@macro\th@definition{\thm@headpunct{}}
\g@addto@macro\th@plain{\thm@headpunct{}}
\makeatother

\usepackage{subfig}

\usepackage[final]{showkeys}

\usepackage{etoolbox}

\usepackage{mathrsfs}

\usepackage[titletoc]{appendix}
\usepackage[doc]{optional}

\usepackage{soul}

\usepackage{colortbl,booktabs,
multirow}
\usepackage{xcolor}

\usepackage{cancel}

\usepackage{empheq}

\definecolor{myblue}{rgb}{.8, .8, 1}

\usepackage{hyperref}
\hypersetup{
colorlinks=true,
linkcolor=blue,
citecolor=blue,
filecolor=magenta,
urlcolor=cyan
}

\usepackage[
open,
openlevel=2,
atend,
numbered
]{bookmark}

\usepackage[capitalize, nameinlink, noabbrev]{cleveref}
\crefname{equation}{}{}
\crefname{chapter}{Chapter}{Chapters}
\crefname{item}{item}{items}
\crefname{figure}{Figure}{Figures}
\crefname{theorem}{Theorem}{Theorems}
\crefname{lemma}{Lemma}{Lemmas}
\crefname{proposition}{Proposition}{Propositions}
\crefname{corollary}{Corollary}{Corollarys}
\crefname{definition}{Definition}{Definitions}
\crefname{fact}{Fact}{Facts}
\crefname{example}{Example}{Examples}
\crefname{algorithm}{Algorithm}{Algorithms}
\crefname{remark}{Remark}{Remarks}
\crefname{note}{Note}{Notes}
\crefname{notation}{Notation}{Notations}
\crefname{case}{Case}{Cases}
\crefname{exercise}{Exercise}{Exercises}
\crefname{question}{Question}{Questions}
\crefname{claim}{Claim}{Claims}
\crefname{enumi}{}{}

\usepackage{pgf}

\usepackage{array}
\usepackage{tabu}

\allowdisplaybreaks

\numberwithin{equation}{section}

\renewcommand\theenumi{(\roman{enumi})}
\renewcommand\theenumii{(\alph{enumii})}

\spnewtheorem*{Proof}{Proof.}{\bf}{\rm}

\def\Ptb{ {\rm \mathbf{P}\mathbf{t}\mathbf{b}}}

\begin{document}

\title{Perturbation Analysis of Error Bounds for Convex Functions on Banach Spaces \thanks{This work was supported by the Natural Science Foundation of Hebei Province (A2024201015), the Excellent Youth Research Innovation Team of Hebei University (QNTD202414) and the Innovation Capacity Enhancement Program-Science and Technology Platform Project, Hebei Province (22567623H) and  benefited from the support of the FMJH Program Gaspard Monge for optimization and operations research and their interactions with data science.}}

\titlerunning{Perturbation Analysis of Error Bounds for Convex Functions on Banach Spaces}

\author{Zhou Wei  \and Michel Th\'era \and Jen-Chih Yao}

\institute{Zhou Wei\at Hebei Key Laboratory of Machine Learning and Computational Intelligence \& College of Mathematics and Information Science, Hebei University, Baoding, 071002, China\\ \email{weizhou@hbu.edu.cn}\\
	Michel Th\'era \at XLIM UMR-CNRS 7252, Universit\'e de Limoges, Limoges, France\ \\ORCID 0000-0001-9022-6406 \\ \email{michel.thera@unilim.fr}\\
	Jen-Chih Yao \at Research Center for Interneural Computing, China Medical University Hospital,
	China Medical University, Taichung, Taiwan \\ and 
	\at Academy of Romanian Scientists, 50044 Bucharest, Romania\\ \email{yaojc@mail.cmu.edu.tw}
}

\dedication{Dedicated to Ralph Tyrrell Rockafellar, with gratitude and admiration,
on the occasion of his 90h birthday}
\date{Received: date / Accepted: date}

\maketitle

\begin{abstract}

This paper focuses on the stability of both local and global error bounds for a proper lower semicontinuous convex function defined on a Banach space. Without relying on any dual space information, we first provide precise estimates of error bound moduli using directional derivatives. For a given proper lower semicontinuous convex function on a Banach space, we prove that the stability of local error bounds under small perturbations is equivalent to the directional derivative at a reference point having a non-zero minimum over the unit sphere. Additionally, the stability of global error bounds is shown to be equivalent to the infimum of the directional derivatives, at all points on the boundary of the solution set, being bounded away from zero over some neighborhood of the unit sphere.


\end{abstract}

\section{Introduction}
This paper focuses on studying the stability of error bounds for lower semicontinuous convex functions defined on a Banach space under data perturbations.
 Given an extended-real-valued function \( f: \mathbb{X} \rightarrow \mathbb{R} \cup \{+\infty\} \) on a Banach space \( \mathbb{X} \), the error bound property of \( f \) is characterized by the following inequality:
\begin{equation}\label{1.1}
c \, {\bf d}(x, \textbf{S}_f) \leq [f(x)]_+,
\end{equation}
where \( \textbf{S}_f := \{x \in \mathbb{X} : f(x) \leq 0\} \) denotes the lower level set of \( f \), \( c \geq 0 \), and \( [f(x)]_+ := \max\{f(x), 0\} \).

 If inequality \eqref{1.1} holds for some \( c > 0 \) and for all \( x \in \mathbb{X} \), then \( f \) is said to have a \textit{global error bound}. For a point \( \bar{x} \in \mathbb{X} \) with \( f(\bar{x}) = 0 \), \( f \) is said to have a \textit{local error bound} at \( \bar{x} \) if inequality \eqref{1.1} holds for some \( c > 0 \) and for all \( x \) sufficiently close to \( \bar{x} \).

 The purpose of this work is to investigate the behavior of error bounds for a given proper lower semicontinuous convex function under perturbations of its data.
 
 It is well-known that Hoffman \cite{28} was a pioneer in the study of error bounds. Specifically, for a given \( m \times n \) matrix \( A \) and an \( m \)-vector \( b \), he showed that the Euclidean distance from a point \( x \) to its nearest point in the polyhedral set \( \{u : Au \leq b\} \) can be bounded above by a constant (dependent only on \( A \)) times the Euclidean norm of the residual error \( \|(Ax-b)_+\| \), where \( (Ax-b)_+ \) denotes the positive part of \( Ax-b \). Hoffman's work has been widely recognized and further developed by many authors, including Robinson~\cite{54}, Mangasarian~\cite{42}, Auslender and Crouzeix \cite{1}, Pang \cite{52}, Lewis and Pang \cite{41}, Klatte and Li \cite{34}, Jourani \cite{jourani}, Abassi and Th\'era \cite{Abassi-Thera1,Abassi-Thera}, Ioffe \cite{ioffe-book}, Mordukhovich \cite{Boris1,Boris2}, and others.

 Error bounds have significant applications in various areas, including sensitivity analysis of linear and integer programs (see \cite{Rob73,Rob77}) and the convergence analysis of descent methods for linearly constrained minimization problems (see \cite{Gul92,HLu,IuD90,TsL92,TsB93}). They are also crucial in addressing the feasibility problem of finding a point in the intersection of finitely many closed convex sets (see \cite{5,6,7,BurDeng02}) and are used in the field of image reconstruction (see \cite{16}).

 This property is closely related to the concepts of weak sharp minima, metric regularity/subregularity of multifunctions, and subtransversality of closed subsets (see \cite{3,BD2,ioffe-JAMS-1,ioffe-JAMS-2,K2,penot-book,54}). For more comprehensive information on the theoretical aspects and practical applications of error bounds, readers can refer to the extensive bibliographies \cite{Aze03,Aze06,7,CKLT,18,DL,22,Kru15,Luke,44,48,penot-book,55,ZTY-JCA,YZ} and the references cited therein.

%
%
%

 This paper focuses on the behavior of error bounds for convex functions when the data of the functions is perturbed, with the goal of exploring the stability of these error bounds under such perturbations. In 1994, Luo and Tseng \cite{LT} conducted a study on the perturbation analysis of a condition number for linear systems. Building on this, Deng \cite{D} explored the perturbation analysis of condition numbers for convex inequality systems. Further research by Az\'e and Corvellec \cite{2} investigated the sensitivity analysis of Hoffman's constants for semi-infinite linear constraint systems. In 2010, Ngai, Kruger, and Th\'era \cite{45} established subdifferential characterizations for the stability of error bounds for semi-infinite convex constraint systems in a Euclidean space. Following this, Kruger, Ngai, and Th\'era \cite{38} studied the stability of error bounds for semi-infinite convex constraint systems in a Banach space, providing dual characterizations related to function perturbations in the system. Kruger, L\'opez, and Th\'era \cite{MP2018} later advanced the subdifferential characterizations of stability presented in \cite{38,45}. In 2022, Wei, Th\'era, and Yao \cite{WTY} characterized the stability of error bounds for convex inequality systems in a Euclidean space using directional derivatives of convex functions. Subsequently, they \cite{WTY2024} extended these results to provide primal characterizations of the stability of error bounds for semi-infinite convex constraint systems in a Banach space.
 
In this paper, we continue the investigation initiated in \cite{WTY2024} on the stability of error bounds for convex functions in a Banach space, with a particular focus on primal characterizations of stability when the functions under consideration are perturbed. We demonstrate that the stability of error bounds for a proper lower semicontinuous convex function is, to some extent, equivalent to solving a class of minimization problems determined by directional derivatives over the unit sphere. Compared to the work presented in \cite{WTY2024}, the main contributions of this paper are twofold:
\begin{itemize}
\item[-] We provide explicit formulae for error bound moduli of a proper lower semicontinuous convex function in terms of directional derivatives and establish a self-contained proof without relying on dual space information (see \cref{th3.1} below). This result is an improvement over \cite[Proposition 3.3]{WTY2024} as it relaxes the continuity assumption made therein.
\item[-] Our stability results lead to primal necessary and sufficient conditions for Hoffman's constants to remain uniformly bounded under perturbations of the problem data (see \cref{th3.5} below). This result improves upon \cite[Theorem 4.2]{WTY2024}.
\end{itemize}

The paper is organized as follows. In Section 2, we present definitions and preliminary results used throughout the paper. Section 3 is devoted to the main results on the stability of local and global error bounds for proper lower semicontinuous convex functions defined on a Banach space. We provide primal characterizations of the stability of error bounds under small perturbations using directional derivatives of component functions (see \cref{th3.2,th3.3,th3.4}). We then apply these stability results to the sensitivity analysis of Hoffman's constant for semi-infinite linear systems and offer a primal equivalent condition for Hoffman's constants to remain uniformly bounded under perturbations (see \cref{th3.5}). Concluding remarks and perspectives are provided in Section 4.
 
\section{Preliminaries}
 
Let $\mathbb{X}$ be a Banach space with topological dual $\mathbb{X}^*$. We denote by $\mathbb{B}_{\mathbb{X}}$ the closed unit ball in $\mathbb{X}$, and for any $x \in \mathbb{X}$ and $\delta > 0$, by $\mathbf{B}(x,\delta)$ the open ball centered at $x$ with radius $\delta$.
 
 Let $D$ be a subset of $\mathbb{X}$. We denote by ${\rm int}(D)$ and ${\rm bdry}(D)$ the interior and boundary of $D$, respectively. The distance from a point $x \in \mathbb{X}$ to the set $D$ is given by
\[
\mathbf{d}(x, D) := \inf_{u \in D} \|x - u\|,
\]
with the convention that $\mathbf{d}(x, D) = +\infty$ whenever $D = \emptyset$.

 In the following, we use the notation $\Gamma_0(\mathbb{X})$ to denote the set of extended-real-valued lower semicontinuous convex functions $\varphi : \mathbb{X} \rightarrow \mathbb{R} \cup \{+\infty\}$, which are assumed to be proper, meaning that ${\rm dom}(\varphi) := \{x \in \mathbb{X} : \varphi(x) < +\infty\}$ is nonempty.
 
 Let $\varphi \in \Gamma_0(\mathbb{X})$ and $\bar{x} \in {\rm dom}(\varphi)$. For any $h \in \mathbb{X}$, we denote by $d^+\varphi(\bar{x}, h)$ the \textit{directional derivative} of $\varphi$ at $\bar{x}$ along the direction $h$, defined by
\begin{equation}\label{2.1}
d^+\varphi(\bar{x}, h) := \lim_{t \rightarrow 0^+} \frac{\varphi(\bar{x} + th) - \varphi(\bar{x})}{t}.
\end{equation}
According to \cite{P}, we know that
\[
t \mapsto \frac{\varphi(\bar{x} + th) - \varphi(\bar{x})}{t}
\]
is nonincreasing as $t \rightarrow 0^+$, and thus
\begin{equation}\label{2.2}
d^+\varphi(\bar{x}, h) = \inf_{t > 0} \frac{\varphi(\bar{x} + th) - \varphi(\bar{x})}{t}.
\end{equation}

 We denote by $\partial \varphi(\bar{x})$ the \textit{subdifferential} of $\varphi$ at $\bar{x}$, defined by
$$
\partial \varphi(\bar{x}) := \{x^* \in \mathbb{X}^* : \langle x^*, x - \bar{x} \rangle \leq \varphi(x) - \varphi(\bar{x}) \ \text{for all } x \in \mathbb{X}\}.
$$
It is known from \cite[Proposition 2.24]{P} that if $\varphi$ is continuous at $\bar{x}$, then $\partial \varphi(\bar{x}) \neq \emptyset$,
\begin{equation}\label{2.3}
\partial \varphi(\bar{x}) = \{x^* \in \mathbb{X}^* : \langle x^*, h \rangle \leq d^+\varphi(\bar{x}, h) \ \text{for all } h \in \mathbb{X}\},
\end{equation}
and
\begin{equation}\label{2.4}
d^+\varphi(\bar{x}, h) = \sup\{ \langle x^*, h \rangle : x^* \in \partial \varphi(\bar{x})\}.
\end{equation}
 
 Given a mapping $\varphi : \mathbb{X} \rightarrow \mathbb{Y}$ from $\mathbb{X}$ to a  \textbf{normed linear space}  $\mathbb{Y}$, we denote by
\[
{\rm Lip}(\varphi ) := \sup_{u, v \in \mathbb{X}, u \neq v} \frac{\|\varphi (u) - \varphi (v)\|_{\mathbb{Y}}}{\|u - v\|_{\mathbb{X}}}
\]
the Lipschitz constant of $\varphi $.

\setcounter{equation}{0}

\section{Stability of Error Bounds for Convex Functions on a Banach Space}

This section is devoted to the stability of local and global error bounds for a given proper lower semicontinuous convex function defined on a Banach space and aims to establish primal characterizations in terms of directional derivatives. We begin with some notations and definitions of error bounds for  convex functions.



Let~$f\in\Gamma_0(X)$. We denote by
\begin{equation}\label{3.1}
	\textbf{S}_f:=\{x\in X: f(x)\leq 0\}
\end{equation}
the {\it lower level set} of $f$. Recall that $f$ is said to have {\it global error bound}, if there exist $c\in(0,+\infty)$ such that
\begin{equation*}
	c {\bf d}(x, \textbf{S}_f)\leq f_+(x)\ \ \forall x\in X,
\end{equation*}
where~$f_+(x):=\max\{f(x), 0\}$. We denote by
\begin{equation}\label{3.2}
	\mathbf{Er} (f):=\inf_{f(x)>0}\frac{f(x)}{{\bf d}(x, \textbf{S}_f)}
\end{equation}
the {\it global error bound modulus} of~$f$. Thus, $f$ has a global error bound if and only if~$ \mathbf{Er} (f)>0$.

Let~$\bar x\in {\rm bdry}(\textbf{S}_f)$. Recall that $f$ is said to have {\it local error bound} at~$\bar x$, if there exist $c,\delta\in(0,+\infty)$ such that
\begin{equation*}
	c {\bf d}(x, \textbf{S}_f)\leq [f(x)]_+\ \ \forall x\in \ball(\bar x,\delta).
\end{equation*}
We denote by
\begin{equation}\label{3.3}
	\mathbf{Er} (f, \bar x):=\liminf_{x\rightarrow\bar x, f(x)>0}\frac{f(x)}{{\bf d}(x, \textbf{S}_f)}
\end{equation}
the {\it local error bound modulus} of~$f$. Thus, ~$f$ has a local error bound at~$\bar x$ ~$ \mathbf{Er} (f, \bar x)>0$.

\medskip

The following theorem, as a key result of this paper, is to give precise estimate on error bound moduli of proper lower semicontinuous convex functions on the Banach space, with  no help of any information from the dual space. The  proof  only relies on the Ekeland variational principle.

\begin{theorem}\label{th3.1}
Let $f\in \Gamma_0(X)$ be such that $\textbf{S}_f$ is nonempty. 
\begin{itemize}
\item[\rm (i)] The following equality on the global error bound modulus holds:
\begin{equation}\label{3.4}
\mathbf{Er}(f)=\inf_{f(x)>0}\left(-\inf_{\|h\|=1}d^+ f(x,h)\right).
\end{equation}
\item[\rm (ii)] Let $\bar x\in{\rm bdry}(S_{f})$. Then the following equality on the local error bound modulus holds:
\begin{equation}\label{3.5}
\mathbf{Er} (f, \bar x)=\liminf_{x\rightarrow\bar x, f(x)>0}\left(-\inf_{\|h\|=1}d^+ f(x,h)\right).
\end{equation}
\end{itemize}
\end{theorem}

{\bf Proof.} (i) Let us define   $$\beta(f):=\inf_{f(x)>0}\left(-\inf_{\|h\|=1}d^+ f(x,h)\right).$$

We first prove that~\eqref{3.4} holds.
Let~$x\in X$ be such that~$f(x)>0$. Take a sequence ~$\{u_n\}$ in $\textbf{S}_f$ such that $\|x-u_n\|\rightarrow {\rm d}(x, \textbf{S}_f)$. By virtue of~\eqref{2.2}, one has 
\begin{eqnarray*}
	\frac{f(x)}{\|x-u_n\|}\leq \frac{f(x)-f(u_n)}{\|x-u_n\|}=-\frac{f(u_n)-f(x)}{\|x-u_n\|}\leq-d^+f\big(x,\frac{u_n-x}{\|x-u_n\|}\big)\leq -\inf_{\|h=1\|}d^+f(x, h).
\end{eqnarray*}
Letting~$n\rightarrow\infty$, one has
$$
\frac{f(x)}{{\bf d}(x, \textbf{S}_f)}\leq-\inf_{\|h=1\|}d^+f(x, h).
$$
This means that
$$
\mathbf{Er} (f)\leq\beta(f).
$$ 

Without loss of generality, we can assume that~$\mathbf{Er} (f)<+\infty$. Suppose on the contrary that~$\mathbf{Er} (f) < \beta(f)$. Then there is~$\tau_0\in\mathbb{R}$ such that~ $\mathbf{Er} (f) < \tau_0 < \beta(f)$. By the definition of~$\mathbf{Er} (f)$, there exists $x_0\in X$ such that 
\begin{equation*}
	f(x_0)>0\ \ {\rm and} \ \ f(x_0)<\tau_0 {\bf d}(x_0, \textbf{S}_f).
\end{equation*}
This implies that
\begin{equation*}
	f_+(x_0)<\inf_{x\in X} f_+(x)+\tau_0 {\bf d}(x_0, \textbf{S}_f).
\end{equation*}
By applying ~Ekeland's variational princinple (cf. \cite{E1974}), there exists ~$y_0\in X$ such that
\begin{equation*}
	\|y_0-x_0\|<{\bf d}(x_0, \textbf{S}_f)\ \ {\rm and} \ \ f_+(y_0)<f_+(x)+\tau_0\|x-y_0\|,\ \ \forall x\in X\backslash\{y_0\}.
\end{equation*}
On the one hand, since     $\|y_0-x_0\|<d(x_0,S_f)$,  this implies that  $y_0\not \in \textbf{S}_f$, that is $f(y_0)=f_+(y_0)>0$. Indeed,  supposing that  $y_0\in S_f$   would imply that 
$\|y_0-x_0\|<d(x_0,S_f)\leq \|y_0-x_0\|$, a contradiction. On the other hand,  for any   ~$h\in X$ with $\|h\|=1$ and ~$t>0$   sufficiently small, one has
$$
\frac{f(y_0+th)-f(y_0)}{t}=\frac{f_+(y_0+th)-f_+(y_0)}{t}>-\tau_0.
$$
Letting $t\rightarrow 0^+$, one yields that $d^+f(y_0, h)\geq -\tau_0$ and consequently
$$
\inf_{\|h\|=1}d^+f(y_0, h)\geq -\tau_0.
$$
Therefore,
$$
-\inf_{\|h\|=1}d^+f(y_0, h)\leq \tau_0<\beta(f)\leq -\inf_{\|h\|=1}d^+f(y_0, h),
$$
which is a contradiction. This means that~$ \mathbf{Er} (f)=\beta(f)$ and so (i) holds.

The same method can be applied to prove \eqref{3.5} and thus we omit the proof of (ii).  The proof is complete.\qed

\vskip 2mm
The following corollary follows from \cref{th3.1} immediately.
\begin{corollary}\label{coro3.1}
Let $f\in \Gamma_0(\mathbb{X})$ and $\bar x\in \mathbb{X}$ such that $f(\bar x)=0$ and $\inf_{\|h\|=1}d^+f(\bar x,h)\not=0$. Then
\begin{equation}\label{3.6}
 \mathbf{Er} (f, \bar x)\geq \left|\inf_{\|h\|=1}d^+f(\bar x,h)\right|.
	\end{equation}
\end{corollary}

{\bf Proof.} If $\inf_{\|h\|=1}d^+f(\bar x, h)>0$, then for any $x\in \mathbb{X}\backslash\{\bar x\}$, one has
\begin{eqnarray*}
	f(x)-f(\bar x)&=&f\Big(\bar x+\|x-\bar x\|\frac{x-\bar x}{\|x-\bar x\|}\Big)-f(\bar x)\\
	&\geq& d^+f\Big(\bar x,\frac{x-\bar x}{\|x-\bar x\|}\Big)\|x-\bar x\|>0.
\end{eqnarray*}
This implies that $S_{f}=\{\bar x\}$ and $\mathbf{Er} (f, \bar x)\geq \inf_{\|h\|=1}d^+f(\bar x,h)$.

If $\inf_{\|h\|=1}d^+f(\bar x, h)<0$, then there exist $t>0$ and $h_0\in \mathbb{X}$ with $\|h_0\|=1$ such that $$f(\bar x+th_0)<0$$ and thus $0\not\in \partial f(\bar x)$ thanks to $f(\bar x)=0$. By virtue of \cite[Theorem 1]{38}, one has  that $f$ admits a local error bound and 
$$\mathbf{Er}(f, \bar x)\geq \mathbf{d}(0, \partial f(\bar x)).$$
Note that for any $x^*\in\partial f(\bar x)$, $h\in \mathbb{X}$ with $\|h\|=1$  and $t>0$, one has
$$
\langle x^*, th\rangle\leq f(\bar x+th)-f(\bar x)
$$
and thus
$$
d^+f(\bar x,h)\geq \langle x^*, h\rangle\geq-\|x^*\|.
$$
This means that
$$
\mathbf{d}(0, \partial f(\bar x))\geq -\inf_{\|h\|=1}d^+ f(\bar x,h)
$$
and consequently 
$$
\mathbf{Er} (f, \bar x)\geq -\inf_{\|h\|=1}d^+f(\bar x,h).
$$
Hence the conclusion holds. The proof is complete.\hfill$\Box$


\begin{remark} It should be noted that inequality in \eqref{3.6} may hold strictly. For example, let $\mathbb{X}:=\mathbb{R}$ and $f(x)\equiv 0$ for all $x\in \mathbb{X}$. Then  one can verify that $$\mathbf{Er} (f, \bar x)=+\infty\ \ {\rm and} \ \ \inf\limits_{\|h\|=1}d^+f(\bar x, h)=0.$$
This means that $\mathbf{Er} (f, \bar x)> \inf_{\|h\|=1}d^+f(\bar x,h)$.
\end{remark}

The following theorem shows that given a proper lower semicontinuous convex function $f$,  the condition $\inf_{\|h\|=1}d^+f(\bar x,h)\not=0$ is essentially equivalent to the stability of local error bounds. For this result, we need the following definition of {\it $\varepsilon$-perturbation of $f$ near $\bar x$} given in \cite[Definition 5]{38}.

\begin{definition}
	Let $f: \mathbb{X}\rightarrow \mathbb{R}\cup\{+\infty\}$, $\bar x\in{\rm dom}(f)$ and $\varepsilon\geq 0$. We say that $g: \mathbb{X}\rightarrow \mathbb{R}\cup\{+\infty\}$ is an $\varepsilon$-perturbation of $f$ near $\bar x$, if
	\begin{equation}\label{3.9}
		\limsup_{x\rightarrow \bar x}\frac{|(f(x)-g(x))-(f(\bar x)-g(\bar x))|}{\|x-\bar x\|}\leq \varepsilon.
	\end{equation}
	In this case, we write $g\in {\rm \mathbf{P}\mathbf{t}\mathbf{b}}(f,\bar x,\varepsilon)$.
\end{definition}
Clearly, if $g\in \Ptb(f,\bar x,\varepsilon)$, then $f\in \Ptb(g,\bar x,\varepsilon)$.

\begin{theorem}\label{th3.2}
Let $f\in \Gamma_0(\mathbb{X})$ and $\bar x\in \mathbb{X}$ with $f(\bar x)=0$. Then the following statements are equivalent:
\begin{itemize}
\item[\rm (i)] $\inf_{\|h\|=1}d^+f(\bar x, h)\not=0$.
\item[\rm (ii)] There exist $c,\varepsilon>0$ such that $\mathbf{Er}(g, \bar x)\geq c$ for all $g\in  \Gamma_0(\mathbb{X})\cap\Ptb(f,\bar x,\varepsilon)$ with $\bar x\in S_g$.

\item[\rm (iii)] There exist $c,\varepsilon>0$ such that $\mathbf{Er}(g, \bar x)\geq c$ for all $g\in  \Gamma_0(\mathbb{X})\cap\Ptb(f,\bar x,\varepsilon)$ defined by 
$$g(\cdot):=f(\cdot)+\varepsilon\langle u^*, \cdot-\bar x\rangle,\ \forall u^*\in\mathbb{B}_{\mathbb{X}^*}.
$$
\end{itemize}
\end{theorem}

{\bf Proof.} We denote $\gamma:=\inf_{\|h\|=1}d^+f(\bar x, h).$

(i) $\Rightarrow$ (ii): Let $\varepsilon,c>0$ be such that $\varepsilon+c<|\gamma|$. Take any $g\in  \Gamma_0(\mathbb{X})\cap\Ptb(f,\bar x,\varepsilon)$ with $\bar x\in S_g$.  \begin{itemize}
\item If $ \gamma>0$, then for any $h\in \mathbb{X}$ with $\|h\|=1$, one has
$$
d^+g(\bar x, h)\geq d^+f(\bar x, h)-\varepsilon, 
$$
and consequently
$$
\inf_{\|h\|=1}d^+g(\bar x, h)\geq\inf_{\|h\|=1}d^+f(\bar x, h)-\varepsilon= \gamma-\varepsilon.
$$
This and \eqref{3.6} imply that $\mathbf{Er}(g, \bar x)\geq \gamma-\varepsilon$.
\item If $\gamma<0$, then for any $h\in \mathbb{X}$, one has
$$
d^+g(\bar x, h)\leq d^+f(\bar x, h)+\varepsilon,
$$
and thus
$$
\inf_{\|h\|=1}d^+g(\bar x, h)\leq\inf_{\|h\|=1}d^+f(\bar x, h)+\varepsilon=\gamma+\varepsilon.
$$
Using \eqref{3.6}, one has $\mathbf{Er}(g, \bar x)\geq -\gamma-\varepsilon>c$. Hence (ii) holds.
\end{itemize}

Note that (ii) $\Rightarrow$ (iii) follows immediately and we next prove (iii) $\Rightarrow$ (i).

Let $\varepsilon>0$. Suppose on the contrary that $\gamma=0$. Then there is a sequence $\{h_k\}\subseteq \mathbb{X}$ such that $\|h_k\|=1$ and
$$
\alpha_k:=d^+f(\bar x,h_k)\rightarrow \gamma=0\ \ {\rm as} \ k\rightarrow\infty.
$$
For sufficiently large $k$, one has $|\alpha_k|<\varepsilon$ and by the Hahn-Banach theorem, there is $h^*_k\in \mathbb{X}^*$ such that $\|h^*_k\|=1$ and $\langle h_k^*, h_k\rangle=\|h_k\|$. 

Let $g(x):=f(x)+\varepsilon \langle h_k^*, x-\bar x\rangle$ for all $x\in \mathbb{X}$. Then $g\in  \Gamma_0(\mathbb{X})\cap\Ptb(f,\bar x,\varepsilon)$. 
Note that $\gamma=0$ and thus for any $x\not=\bar x$, one has $f(x)\geq f(\bar x)$. By \eqref{2.1}, there exists $\delta_k\rightarrow 0^+$ such that
\begin{equation}\label{3.10}
  f(\bar x+\delta_kh_k)<f(\bar x)+(\varepsilon+\alpha_k)\delta_k=\inf_{x\in \mathbb{X}}f(x)+(\varepsilon+\alpha_k)\delta_k.
\end{equation}
By virtue of  Ekeland's variational principle (cf. \cite{E1974}), there is $z_k\in \mathbb{X}$ such that
$$
\|z_k-(\bar x+\delta_kh_k)\|<\frac{\delta_k}{2}, f(z_k)\leq f(\bar x+\delta_kh_k)
$$
and
\begin{equation}\label{3.11}
 f(x)+2(\varepsilon+\alpha_k)\|x-z_k\|>f(z_k),\ \forall x: x\not=z_k.
\end{equation}
This implies that $z_k\rightarrow \bar x$, $g(\bar x)=f(\bar x)=0$ and
\begin{eqnarray*}
g(z_k)&=&f(z_k)+\varepsilon\langle h_k^*, z_k-\bar x \rangle\\
&\geq& f(\bar x)+\varepsilon\langle h_k^*, z_k-\bar x \rangle\\
&=&\varepsilon\langle h_k^*, z_k-\bar x-\delta_kh_k \rangle+\varepsilon\delta_k\\
&>&\varepsilon\delta_k-\frac{1}{2}\varepsilon\delta_k=\frac{1}{2}\varepsilon\delta_k>0.
\end{eqnarray*}
We claim that
\begin{equation}\label{3.12}
  \inf_{\|h\|=1}d^+g(z_k,h)< 0,
\end{equation}
Indeed, by  reductio ad absurdum, if  $\inf_{\|h\|=1}d^+g(z_k,h)\geq 0$,  thus $g(z_k)=\inf_{x\in \mathbb{X}}g(x)$, which contradicts $g(\bar x)=0$.

For any $h\in \mathbb{X}$, $\|h\|=1$ and any $t>0$, by \eqref{3.11}, one can get
\begin{eqnarray*}
\frac{g(z_k+th)-g(z_k)}{t}=\frac{f(z_k+th)-f(z_k)}{t}+\varepsilon\langle h_k^*, h\rangle\geq-2(\varepsilon+\alpha_k)\|h\|-\varepsilon=-5\varepsilon.
\end{eqnarray*}
Thus,
$$
\inf_{\|h\|=1}d^+g(z_k,h)\geq-5\varepsilon.
$$
By \eqref{3.12} and \cref{th3.1}, one has  $\mathbf{Er}(g, \bar x)\leq5\varepsilon\rightarrow 0^+$, which contradicts (iii), as $\varepsilon$ is arbitrary. The proof is complete. \hfill$\Box$

\medskip

The following corollary holds from \cref{th3.2} immediately.

\begin{corollary}
Let $f\in \Gamma_0(\mathbb{X})$ and $\bar x\in \mathbb{X}$ with $f(\bar x)=0$. Then 
\begin{itemize}
	\item [\rm(i)] $\inf_{\|h\|=1}d^+f(\bar x,h)\not=0$ holds if and only if there exists $\varepsilon>0$ such that
	$$
	\inf\big\{\mathbf{Er}(g, \bar x): g\in  \Gamma_0(\mathbb{X})\cap\Ptb(f,\bar x,\varepsilon)\ {\it with} \ \bar x\in S_g \big\}>0.
	$$
	\item [\rm(ii)] $\inf_{\|h\|=1}d^+f(\bar x,h)=0$ holds if and only if for any $\varepsilon>0$, one has
	$$
	\inf\big\{\mathbf{Er}(g, \bar x): g\in \Gamma_0(\mathbb{X})\cap \Ptb(f,\bar x,\varepsilon)\ {\it with} \ \bar x\in S_g \big\}=0.
	$$
\end{itemize}

\end{corollary}

Next, we study the stability of global error bounds for lower semicontinuous convex functions. The following theorem give  primal criterion for the stability of global error bounds.


\begin{theorem}\label{th3.3}
Let $f\in \Gamma_0(\mathbb{X})$ be such that ${\rm bdry}(\textbf{S}_f)\subseteq f^{-1}(0)$. Consider the following statements:
\begin{itemize}

\item[\rm (i)] there exist $c,\varepsilon\in (0, +\infty)$ such that  $\mathbf{Er}(g)\geq c$ holds for all $g\in \Gamma_0(\mathbb{X})$ satisfying
    \begin{equation}\label{3-10}
      \textbf{S}_f\subseteq S_g\  \ {\it and} \ \ {\rm Lip}(f-g)<\varepsilon;
    \end{equation}
\item[\rm (ii)] there exist $c,\varepsilon\in (0, +\infty)$ such that $\mathbf{Er}(g)\geq c$ holds for all $g\in  \Gamma_0(\mathbb{X})$ satisfying
    \begin{equation}\label{3-11}
      {\rm bdry}(\textbf{S}_f)\cap g^{-1}(0)\not=\emptyset\ \  {\it and}  \  \ {\rm Lip}(f-g)<\varepsilon.
    \end{equation}
\end{itemize}
Then {\rm (i)} (resp. {\rm (ii)}) holds if (resp. only if) there is $\tau\in (0, +\infty)$ such that
\begin{equation}\label{3-9}
	\inf\left\{\Big|\inf_{\|h\|=1}d^+f(\bar x,h)\Big|: \bar x\in{\rm bdry}(\textbf{S}_f)\right\}\geq\tau.
\end{equation}
\end{theorem}

{\bf Proof.} {\it If part}: Suppose that there is $\tau\in (0, +\infty)$ such that \eqref{3-9} holds. 
First  if  there is some $\bar x\in{\rm bdry}(\textbf{S}_f)$ such that $\inf_{\|h\|=1}d^+f(\bar x, h)>0$, using the proof of  \cref{coro3.1}, the lower level set $\textbf{S}_f$ is a  singleton and thus $\textbf{S}_f=\{\bar x\}$. This means that the conclusion in (ii) follows by using \cref{th3.2} as well as its proof.

Next, we assume that $\inf_{\|h\|=1}d^+f(\bar x, h)\leq 0$ for all $\bar x\in{\rm bdry}(\textbf{S}_f)$. We claim that
\begin{equation}\label{3-12}
  \tau\mathbf{d}(x, \textbf{S}_f)\leq [f(x)]_+,\ \ \forall x\in \mathbb{X}.
\end{equation}
Let $x\in\mathbb{X}\backslash \textbf{S}_f$ and $\nu>0$ sufficiently small. Take $x_{\nu}\in {\rm bdry}(S_{f})$ such that
\begin{equation}\label{3.017a}
  \|x-x_{\nu}\|<(1+\nu)\mathbf{d}(x, \textbf{S}_f).
\end{equation}
By virtue of \cref{th3.1}, one has
$$
\mathbf{Er}(f,x_{\nu})=\liminf_{u\rightarrow x_{\nu}, u\not\in \textbf{S}_f}\Big(-\inf_{\|h\|=1}d^+f(u, h)\Big)\geq-\inf_{\|h\|=1}d^+f(x_{\nu}, h)\geq\tau.
$$

Using the definition of $\mathbf{Er}(f,x_{\nu})$, for $\nu>0$, there is $r>0$ such that for any $u\in {\bf B}(x_{\nu},r)\backslash\textbf{S}_f$, one has
$$
\frac{f(u)}{\mathbf{d}(u, \textbf{S}_f)}>\mathbf{Er}(f,x_{\nu})-\nu\geq \tau-\nu.
$$
Note that $x\not\in \textbf{S}_f$ and $x_{\nu}\in{\rm bdry}(\textbf{S}_f)$, and thus
$$
\{\lambda x_{\nu}+(1-\lambda)x: \lambda\in (0,1)\}\cap \textbf{S}_f=\emptyset.
$$
Then there is $\lambda\in (0,1)$ such that $\lambda x_{\nu}+(1-\lambda)x\in B(x_{\nu},r)\backslash\textbf{S}_f$

and therefore
\begin{eqnarray*}
\tau\leq\frac{f(\lambda x_{\nu}+(1-\lambda)x)}{\mathbf{d}(\lambda x_{\nu}+(1-\lambda)x, \textbf{S}_f)}+\nu\leq\frac{\lambda f(x_{\nu})+(1-\lambda)f(x)}{\mathbf{d}(x, \textbf{S}_f)-\|\lambda x_{\nu}+(1-\lambda)x-x\|}+\nu,
\end{eqnarray*}
and it follows from \eqref{3.017a} that

\begin{eqnarray*}
\tau\leq\frac{(1-\lambda)f(x)}{(1-\lambda(1+\nu))\mathbf{d}(x, \textbf{S}_f)}  \rightarrow \frac{f(x)}{\mathbf{d}(x, \textbf{S}_f)}\ {\rm as\ }\nu\rightarrow 0^+.
\end{eqnarray*}
This means that \eqref{3-12} holds.

Let $\varepsilon\in (0,\tau)$ and $g\in \Gamma_0(\mathbb{X})$ satisfying \eqref{3-10}. Take any $x\in \mathbb{X}$ such that $g(x)>0$. Then $f(x)>0$ (thanks to $\textbf{S}_f\subseteq S_g$).

We claim that
\begin{equation}\label{3.13}
  \inf_{\|h\|=1}d^+f(x, h)\leq -\tau.
\end{equation}
Granting this, it follows from ${\rm Lip}(f-g)<\varepsilon$ in \eqref{3-10} that
$$
\inf_{\|h\|=1}d^+g(x,h)\leq\inf_{\|h\|=1}d^+f(x,h)+\varepsilon\leq -(\tau-\varepsilon).
$$
This and \cref{th3.1} imply that $\mathbf{Er}(g)\geq \tau-\varepsilon$.

It remains to prove \eqref{3.13}. For any  $n\in\mathbb{N}$, we can take $z_n\in {\rm bdry}(\textbf{S}_f)$ such that
$$
\|x-z_n\|<(1+\frac{1}{n})\mathbf{d}(x,\textbf{S}_f).
$$
Then \eqref{3-12} implies that
$$
f(x)\geq \tau \mathbf{d}(x, \textbf{S}_f)>\frac{n}{n+1}\tau \|x-z_n\|.
$$
For any $t\in (0,1)$, one has
$$
f(x+t(z_n-x))\leq tf(z_n)+(1-t)f(x),
$$
and thus
$$
\frac{f(x+t(z_n-x))-f(x)}{t}\leq -f(x)\leq-\frac{n}{n+1}\tau\|x-z_n\|.
$$
This implies that
$$
\inf_{\|h\|=1}d^+f(x, h)\leq d^+f\Big(x,\frac{z_n-x}{\|x-z_n\|}\Big)\leq-\frac{n}{n+1}\tau.
$$
Hence
$$
\inf_{\|h\|=1}d^+f(x, h)\leq-\frac{n}{n+1}\tau.
$$
Letting $n\rightarrow\infty$, one has that \eqref{3.13} holds.

{\it Only if part}: we assume that (ii) holds. Suppose on  the contrary that there is a sequence $\{x_k\}\subseteq {\rm bdry}(\textbf{S}_f)$ such that
$$
\alpha_k:=\inf_{\|h\|=1}d^+f(x_k, h)\rightarrow 0^- \ \ {\rm as} \ k\rightarrow \infty.
$$
Let $\varepsilon>0$ and take   a  sufficiently large $k$ such that
\begin{equation}\label{3.14}
  \frac{3}{2}\alpha_k+\frac{\varepsilon}{2}>0.
\end{equation}
Note that for any $x\not=x_k$, one has
$$
\frac{f(x)-f(x_k)}{\|x-x_k\|}=\frac{f\big(x_k+\|x-x_k\|\cdot\frac{x-x_k}{\|x-x_k\|}\big)-f(x_k)}{\|x-x_k\|}\geq d^+f\Big(x_k,\frac{x-x_k}{\|x-x_k\|}\Big)\geq\alpha_k,
$$
and consequently
\begin{equation}\label{3.15}
  f(x)-\alpha_k\|x-x_k\|\geq f(x_k),\ \ \forall x\in \mathbb{X}.
\end{equation}
By the definition of $\alpha_k$, there exists $h_k\in \mathbb{X}$ such that
\begin{equation}\label{3.16}
\|h_k\|=1 \ \ {\rm and} \ \  d^+f(x_k,h_k)<\alpha_k+\frac{\varepsilon}{2}.
\end{equation}
Then one can choose $r_k\rightarrow 0^+$ (as $k\rightarrow \infty$) such that
\begin{equation}\label{3.17}
  f(x_k+r_kh_k)<f(x_k)+(\alpha_k+\varepsilon)r_k.
\end{equation}
This and \eqref{3.15} imply that
$$
  f(x_k+r_kh_k)-\alpha_k\|x_k+r_kh_k-x_k\|<\inf_{x\in X}(f(x)-\alpha_k\|x-x_k\|)+\varepsilon r_k.
$$
Applying  Ekeland's variational principle (cf. \cite{E1974}), there exists $y_k\in \mathbb{X}$ such that
\begin{equation}\label{3.18}
  \|y_k-(x_k+r_kh_k)\|<\frac{r_k}{2}, f(y_k)-\alpha_k\|y_k-x_k\|\leq  f(x_k+r_kh_k)-\alpha_k r_k,
\end{equation}
and
\begin{equation}\label{3.19}
  f(x)-\alpha_k\|x-x_k\|+2\varepsilon\|x-y_k\|>f(y_k)-\alpha_k\|y_k-x_k\|,\ \ \forall x: x\not=y_k.
\end{equation}
Then
$$
\|y_k-x_k\|> r_k-\frac{r_k}{2}=\frac{r_k}{2}\ \ {\rm and}\ \ \|y_k-x_k\|< r_k+\frac{r_k}{2}=\frac{3}{2}r_k,
$$
and thus $y_k\not=x_k$.

For any $h_k$, by  the Hahn-Banach theorem, there exists $h^*_k\in \mathbb{X}^*$ such that
$$
\|h^*_k\|=1\ \ {\rm and} \ \ \langle h_k^*, h_k\rangle=\|h_k\|.
$$
We consider the following function $g$:
$$
g(x):=f(x)+\varepsilon\langle h_k^*, x-x_k\rangle\ \ \forall x\in \mathbb{X}.
$$
Then $g\in  \Gamma_0(\mathbb{X})$. By virtue of \eqref{3.14}, \eqref{3.15}, \eqref{3.16} and \eqref{3.19}, one has
\begin{eqnarray*}
g(y_k)=f(y_k)+\varepsilon\langle h_k^*, y_k-x_k\rangle&=&f(y_k)+\varepsilon\langle h_k^*, y_k-(x_k+r_kh_k)\rangle+\varepsilon r_k\\
&\geq&\alpha_k\|y_k-x_k\|-\varepsilon\|y_k-(x_k+r_kh_k)\|+\varepsilon r_k\\
&\geq&\alpha_k\cdot\frac{3}{2}r_k+\frac{\varepsilon}{2} r_k>0.
\end{eqnarray*}
\vskip 2mm
We consider two cases:
\begin{itemize}
\item If $\inf_{\|h\|=1}d^+g(y_k, h)\geq 0$, then for any $x\not=y_k$, one has
$$
g(x)-g(y_k)\geq d^+g\big(y_k,\frac{x-y_k}{\|x-y_k\|}\big)\|x-y_k\|\geq \inf_{\|h\|=1}d^+g(y_k, h)\|x-y_k\|\geq 0.
$$
This together with $g(y_k)>0$ implies that  $S_{g}=\emptyset$ and thus $\mathbf{Er}(g)=0$, which contradicts (ii).
\item If  $\inf_{\|h\|=1}d^+g(y_k, h)<0$. For any $h\in \mathbb{X}$ with $\|h\|=1$ and any  $t>0$, by \eqref{3.19}, one has
\begin{eqnarray*}
\frac{g(y_k+th)-g(y_k)}{t}&=&\frac{f(y_k+th)-f(y_k)}{t}+\varepsilon\langle h_k^*, h\rangle\\
&\geq&\frac{1}{t}\big(\alpha_k\|y_k+th-x_k\|-\alpha_k\|y_k-x_k\|-2\varepsilon \|y_k+th-y_k\|\big)+\varepsilon\langle h_k^*, h\rangle\\
&\geq&\alpha_k-2\varepsilon-\varepsilon,
\end{eqnarray*}
and consequently
$$
0>\inf_{\|h\|=1}d^+g(y_k, h)\geq\alpha_k-2\varepsilon-\varepsilon\geq-4\varepsilon.
$$
Then \cref{th3.1} gives that $\mathbf{Er}(g)\leq 4\varepsilon$, which contradicts (ii), as $\varepsilon$ is arbitrary. 
\end{itemize}
The proof is complete. \hfill$\Box$



\begin{remark}\label{re3-2} 
It is noted that the validity of condition \eqref{3-9} (for some $\tau>0$)  may not be sufficient for the stability of the global error bound for $f$  as in (ii). For example, let $f(x):=e^x-1$ for all $x\in \mathbb{R} $. Then $f(x)=e^x-1$, $\textbf{S}_f=(-\infty, 0]$, ${\rm bdry}(\textbf{S}_f)=\{0\}$ and $|\inf_{|h|=1}d^+f(0, h)|=1>0$. Let $\varepsilon\in (0,+\infty)$, and consider $g_{\varepsilon}(x):=f(x)-\varepsilon x, \forall x\in\mathbb{R}$. Then 
	for any $x<\bar x$, one has
$$
\frac{g_{\varepsilon}(x)}{\mathbf{d}(x, S_{g_{\varepsilon}})}=\frac{e^x-1-\varepsilon x}{\bar x-x}\rightarrow \varepsilon \ \ {\rm as} \ x\rightarrow -\infty.
$$
This implies that $\mathbf{Er}(g_{\varepsilon})\leq 2\varepsilon,$
and consequently the stability of global error bound for $f$  as said in (ii) does not hold (as $\varepsilon>0$ is arbitrary). This gives rise to the question:  {\it 
	 whether  there  is some type of stability of global error bounds that can be characterized by condition \eqref{3-9}? } We have no answer to this question.
\end{remark}

The following theorem shows that when adding  to \eqref{3-9} a mild  qualification assumption  makes it possible to characterize the stability of the global error bound for $f$  as said in (ii) of \cref{th3.4}.

\vskip 2mm

\begin{theorem}\label{th3.4}
Let $f\in \Gamma_0(\mathbb{X})$ be such that ${\rm bdry}(\textbf{S}_f)\subseteq f^{-1}(0)$.  Then the following statements are equivalent:
\begin{itemize}
\item[\rm (i)] There exists $\tau\in (0, +\infty)$ such that \eqref{3-9} holds and
\begin{equation}\label{3.20}
\liminf_{k\rightarrow\infty}\Big|\inf_{\|h\|=1}d^+f(z_k,h)\Big|\geq\tau, \ \forall \{(z_k,x_k)\}\subseteq {\rm int}(\textbf{S}_f)\times{\rm bdry}(\textbf{S}_f)\ {\it with} \ \lim\limits_{k\rightarrow\infty}\frac{f(z_k)-f(x_k)}{\|z_k-x_k\|}=0;
\end{equation}
\item[\rm (ii)] there exist $c,\varepsilon\in (0, +\infty)$ such that $\mathbf{Er}(g)\geq c$ holds for all $g\in  \Gamma(\mathbb{X})$ satisfying \eqref{3-11};
\item[\rm (iii)] there exist $c,\varepsilon>0$ such that $\mathbf{Er}(g_{u^*, \varepsilon})\geq c$ holds for all $g\in  \Gamma(\mathbb{X})$ defined by 
$$g_{u^*, \varepsilon}(\cdot):=f(\cdot)+\varepsilon\langle u^*, \cdot-\bar x\rangle, \forall (\bar x,u^*)\in{\rm bdry}(\textbf{S}_f) \times \mathbb{B}_{\mathbb{X}^*}.
$$
\end{itemize}
\end{theorem}

{\bf Proof.} (i)$\Rightarrow$(ii): If there is $\bar x\in{\rm bdry}(\textbf{S}_f)$ such that $\inf_{\|h\|=1}d^+f(\bar x, h)>0$, the conclusion follows by \cref{th3.1} as well as the proof of  \cref{th3.2}.

Next, we consider $\inf_{\|h\|=1}d^+f(\bar x, h)\leq 0$ for all $\bar x\in{\rm bdry}(\textbf{S}_f)$. We first prove the following claim:
\begin{claim}
  There exists $\varepsilon_0>0$ such that for any $x_0\in{\rm bdry}(\textbf{S}_f)$, one has
\begin{equation}\label{3.21}
\inf\left\{\Big|\inf_{\|h\|=1}d^+f(z_0,h)\Big|: z_0\in \mathbb{X}, f(z_0)\geq-\varepsilon_0\|z_0-x_0\|\right\}\geq\tau.
\end{equation}
\end{claim}
\textbf{Proof of the claim}

Suppose on the contrary that there exist $\varepsilon_k\rightarrow 0^+$, $x_k\in{\rm bdry}(\textbf{S}_f)$ and $z_k\in\mathbb{X}$ such that
\begin{equation}\label{3.22}
  f(z_k)\geq-\varepsilon_k\|z_k-x_k\|\ \ {\rm and} \ \ \Big|\inf_{\|h\|=1}d^+f(z_k,h)\Big|<\tau,\ \ \forall  k.
\end{equation}
Then $f(z_k)\leq 0$  for all $k$, since otherwise, using the proof of \eqref{3.13}, one get $$\big|\inf_{\|h\|=1}d^+f(z_k,h)\big|>\tau,
 $$
which contradicts \eqref{3.22}.

By virtue of \eqref{3-9}, one has $z_k\in \textbf{S}_f\backslash {\rm bdry}(\textbf{S}_f)$ and then \eqref{3.22} implies
$$
0\geq\frac{f(z_k)-f(x_k)}{\|z_k-x_k\|}=\frac{f(z_k)}{\|z_k-x_k\|}\geq-\varepsilon_k.
$$
By \eqref{3.20}, one gets
$$
\liminf_{k\rightarrow\infty}\Big|\inf_{\|h\|=1}d^+f(z_k,h)\Big|>\tau,
$$
which contradicts \eqref{3.22}. Hence the claim holds. \qed
\vskip 2mm
Let $\varepsilon>0$ be such that $\varepsilon<\min\{\varepsilon_0, \tau\}$ and $g\in  \Gamma_0(\mathbb{X})$ be such that \eqref{3-11} holds. Take any $\bar x\in{\rm bdry}(\textbf{S}_f)\cap g^{-1}(0)$. Then for any $x\in \mathbb{X}$  with $g(x)>0$, one has
$$
f(x)\geq g(x)+(f(\bar x)-g(\bar x))-\varepsilon\|x-\bar x\|>-\varepsilon\|x-\bar x\|.
$$
By \eqref{3.21}, one has
$$
\inf_{\|h\|=1}d^+f(x,h)<-\tau,
$$
and consequently
$$
\inf_{\|h\|=1}d^+g(x,h)<\inf_{\|h\|=1}d^+f(x,h)+\varepsilon<-(\tau-\varepsilon).
$$
This and \cref{th3.1} imply that $\mathbf{Er} g\geq \tau-\varepsilon$.

Note that  (ii)$\;\Rightarrow\;$(iii) follows immediately and it remains to prove (iii)$\;\Rightarrow\;$(i).

Suppose on the contrary that the conclusions in (i) does not hold. Based on \cref{th3.4} and its proof, it suffices to assume that \eqref{3.20} does not hold for any $\tau>0$. Then, there exists $(z_k,x_k)\in {\rm int}(\textbf{S}_f)\times {\rm bdry}(\textbf{S}_f)$ such that
\begin{equation}\label{3.23}
  \lim_{k\rightarrow\infty}\frac{f(z_k)-f(x_k)}{\|z_k-x_k\|}=0\ \ {\rm and} \ \ \alpha_k:=\inf_{\|h\|=1}d^+f(z_k,h)\rightarrow 0^- \ \ {\rm as}\ k\rightarrow\infty.
\end{equation}
Let $\varepsilon>0$. For any $k$, by the Hahn-Banach theorem,  there is $h_k^*\in \mathbb{X}^*$ such that
$$
\|h_k^*\|=1\ \ {\rm and} \ \ \langle h_k^*, z_k-x_k\rangle=\|z_k-x_k\|.
$$
Then when $k$ is sufficiently large, one has
\begin{equation}\label{3.24}
   \alpha_k+\varepsilon>0\ \  {\rm and} \ \ \frac{f(z_k)-f(x_k)}{\|z_k-x_k\|}+\varepsilon\Big\langle h_k^*, \frac{z_k-x_k}{\|z_k-x_k\|}\Big\rangle>0.
\end{equation}
Consider the function $g_{k,\varepsilon}$ defined as follows:
$$
g_{k,\varepsilon}(x):=f(x)+\varepsilon\langle h_k^*, x-x_k\rangle\ \ \forall x\in \mathbb{X}.
$$
Then $g_{k,\varepsilon}\in   \Gamma_0(\mathbb{X})$ and it follows from \eqref{3.24} that $$g_{k,\varepsilon}(z_k)=f(z_k)+\varepsilon\langle h_k^*, z_k-x_k\rangle>0.$$
Thus
$$
0>\inf_{\|h\|=1}d^+g_{k,\varepsilon}(z_k,h)\geq\inf_{\|h\|=1}d^+f(z_k,h)-\varepsilon
=\alpha_k-\varepsilon>-2\varepsilon.
$$
Then \cref{th3.1} gives $\mathbf{Er} g_{k,\varepsilon}\leq 2\varepsilon$, which contradicts (iii), as $\varepsilon>0$ is arbitrary. The proof is complete.\qed

\begin{remark} It should be noted that  the stability of the global error bound may not hold if condition \eqref{3.20} is violated. Consider again the example given in \cref{re3-2}. Then $f(x)=e^x-1$, $\textbf{S}_f=(-\infty, 0]$, ${\rm bdry}(\textbf{S}_f)=\{0\}$ and the stability of the global error bound in \cref{th3.5} is not satisfied. Further, for any $z_k\rightarrow -\infty$, one can verify that
$$
\Big|\inf_{|h|=1}d^+f(z_k,h)\Big|=e^{z_k}\rightarrow 0 \ \ {\rm as} \ k\rightarrow \infty,
$$
which implies that \eqref{3.20} does not hold. \qed
\end{remark}

\section{Application to the sensibility analysis of  Hoffman's constants for semi-infinite linear systems}

As an application of \cref{th3.4}, we study sensitivity analysis of Hoffman's constants for semi-infinite linear systems in a Banach space and aim  to establish a primal characterization for the Hoffman's constants to be uniformly bounded under perturbations on the problem data.

Let $T$ be a compact, possibly infinite metric space and $a^*:T\rightarrow \mathbb{X}^*, b:T\rightarrow \mathbb{R}$ be continuous functions on $T$. We consider now semi-infinite linear systems in $\mathbb{X}$ defined by
\begin{equation}\label{4-1a}
  \langle a^*(t), x\rangle \leq \ball(t),\ \ {\rm for\ all} \ t\in T.
\end{equation}
We denote by $\mathcal{S}_{a^*,b}$ the set of solutions to system \eqref{4-1a}. We use the following notations:
$$
\begin{aligned}
&f_{a^*,b}(x):=\max_{t\in T}( \langle a^*(t), x\rangle -\ball(t)),\\
&J_{a^*,b}(x):=\{t\in T: \langle a^*(t), x\rangle -\ball(t))=f_{a^*,b}(x)\} \ \  {\rm for \ each } \ x\in \mathbb{X}.
\end{aligned}
$$
It is easy to verify that $J_{a^*,b}(x)$ is a compact subset of $T$ for each $x\in\mathbb{X}$.

Recall that $\mathcal{S}_{a^*,b}$ admits a global error bound, if there exists $\sigma>0$ such that
\begin{equation}\label{4-2a}
 \sigma \mathbf{d}(x, \mathcal{S}_{a^*,b})\leq [f_{a^*,b}(x)]_+, \  {\rm for\ all} \ x\in \mathbb{X}.
\end{equation}
The Hoffman constant of the semi-infinite linear system \eqref{4-1a}, denoted by $\sigma(a^*,b)$, is given by
\begin{equation}\label{4-3a}
\sigma(a^*,b):=\sup\{\sigma>0:  \text{such that } \eqref{4-2a}\ {\rm holds}\}.
\end{equation}
One can verify that
$$
\sigma(a^*,b)=\inf_{x\not\in \mathcal{S}_{a^*,b}}\frac{f_{a^*,b}(x)}{\mathbf{d}(x, \mathcal{S}_{a^*,b})}.
$$

Note that $f_{a^*,b}(\cdot)$ is continuous and thus one can verify that 
\begin{equation}\label{3.31}
	{\rm bdry}(\mathcal{S}_{a^*,b})=f_{a^*,b}^{-1}(0).
\end{equation}


\vskip 2mm

The following theorem gives a primal equivalent criteria on the Hoffman's constant for system \eqref{4-1a} under perturbations on the problem data.

\begin{theorem}\label{th3.5}
Let $\mathcal{J}:=\{J_{a^*,b}(x): f_{a^*,b}(x)=0\}$. Then there exists $\tau>0$ such that
\begin{equation}\label{3.32}
\inf\left\{\Big|\inf_{\|h\|=1}\sup_{t\in J}\langle a^*(t), h\rangle\Big|: J\in \mathcal{J} \right\}\geq \tau
\end{equation}
if and only if there exist $c,\varepsilon>0$ such that $\sigma(\widetilde a^*,\widetilde b)\geq c$ holds for all $(\widetilde a^*,\widetilde b )$ defined by 
    \begin{equation}\label{4-7a}
      \widetilde a^*(t):=a^*(t)+\varepsilon \widetilde u^*\ \ {\it and} \ \ \widetilde \ball(t):=\ball(t)+\varepsilon\langle \widetilde u^*,\widetilde x\rangle, \  \ \forall (\widetilde x, \widetilde u^*)\in{\rm bdry}(\mathcal{S}_{a^*,b})\times  \mathbb{B}_{\mathbb{X}^*}.
          \end{equation}
Moreover, the following inequality on Hoffman's constant $\sigma(a^*,b)$ holds:
\begin{equation}\label{4-8a}
\sigma(a^*,b)\geq\inf\left\{\Big|\inf_{\|h\|=1}\sup_{t\in J}\langle a^*(t), h\rangle\Big|: J\in \mathcal{J}\right\}.
\end{equation}

\end{theorem}

{\bf Proof.} By virtue of \cref{th3.4}, we only need to show that $f_{a^*,b}$ and $\tau$ satisfy assumption \eqref{3.20} in \cref{th3.4}.

Indeed, for each  $x\in {\rm bdry}(\mathcal{S}_{a^*,b})$, one has that $J_{a^*,b}(x)\in \mathcal{J}$ and it follows from \cite[Theorem 4.2.3]{SIAM13}(or \cite[Proposition 4.5.2]{S2007}) that
$$
\inf_{\|h\|=1}d^+f_{a^*,b}(x,h)=\inf_{\|h\|=1}\sup_{t\in J_{a^*,b}(x)}\langle a^*(t), h\rangle.
$$
This and \eqref{3.32} imply that
$$
\inf_{\|h\|=1}d^+f_{a^*,b}(x,h)\geq\tau$$
and thus $f$ and $\tau$ satisfy \eqref{3.20} as said in \cref{th3.4}. The proof is complete.\qed

\vskip 2mm

As an application of \cref{th3.5}, we give one examples to verify the uniform bound of the Hoffman constant for the given linear system.

\begin{example}
Consider the following linear system \eqref{4-1a}, where $\mathbb{X}:=\mathbb{R}^2$, $T:=\{1,2,3\}$, and $a^*:T\rightarrow \mathbb{R}^2, b:T\rightarrow \mathbb{R}$ are defined as follows:
$$
a^*(1):=(1,1), a^*(2):=(-2,1), a^*(3):=(1,-2), \ball(1):=1, \ball(2):=2, \ball(3):=2.
$$
Then one can verify that
$$
\mathcal{J}=\left\{\{1\},\{2\}, \{3\}, \{1,2\}, \{2,3\}, \{1,3\}\right\}.
$$
Solve the following six min-max optimization problems:
\begin{equation}
OP(J) \ \ \ \
 \min_{\|h\|=1} \max_{t\in J_k}\langle a^*(t),h\rangle, J\in \mathcal{J}.
\end{equation}
Then one has
$$
\min\left\{\Big|\min_{\|h\|=1}\max_{t\in J}\langle a^*(t), h\rangle\Big|: J\in \mathcal{J}\right\}=\frac{\sqrt{2}}{2}>0,
$$
and then \cref{th3.5} implies that the Hoffman constant ($\sigma(\cdot)$) for system \eqref{4-1a} is uniformly bounded under perturbations on the problem data $(a^*,b)$.

\end{example}

Finally, we  show  that the uniform bound of the Hoffman constant may not be satisfied if the condition \eqref{3.32} is violated in  \cref{th3.5}.

\begin{example}
Consider the following linear system \eqref{4-1a}, where  $\mathbb{X}:=\mathbb{R}^2$, $T:=\{1,2\}$, and $a^*:T\rightarrow \mathbb{R}^2, b:T\rightarrow \mathbb{R}$ are defined as follows:
	$$
	a^*(1):=(1,1), a^*(2):=(-1,-1), \ball(1)=\ball(2):=0.
	$$
	Then $\mathcal{S}(a^*,b)=\{x=(x_1, x_2)\in \mathbb{R}^2: x_1+x_2=0\}$. Note that
	\begin{equation*}
	J_{a^*,b}(0)=T,	\min_{\|h\|=1} \max_{t\in T}\langle a^*(t),h\rangle=0
	\end{equation*}
	and thus condition \eqref{3.32} is violated.
	
	 We next show that the Hoffman constant for system \eqref{4-1a} is not uniformly bounded under perturbations on the problem data $(a^*,b)$.
	
	Indeed, for any $\varepsilon>0$,   select $u_{\varepsilon}^*:=(0,\varepsilon)$, $x_{\varepsilon}:=(0,0)$ and define $a_{\varepsilon}^*(t):=a^*(t)+u_{\varepsilon}^*,b_{\varepsilon}(t):=\ball(t)$ for $t=1,2$. Then $f_{a_{\varepsilon}^*,b_{\varepsilon}}(x):=\max\{x_1+x_2+\varepsilon x_2, -x_1-x_2+\varepsilon x_2\}$ for all $x=(x_1, x_2)\in \mathbb{R}^2$ and
	$$
	\mathcal{S}(a_{\varepsilon}^*,b_{\varepsilon})=\{x=(x_1, x_2)\in \mathbb{R}^2: x_1+(1+\varepsilon)x_2\leq 0, x_1+(1-\varepsilon)x_2\geq 0\}.
	$$
	Take $u_{\varepsilon}:=(-\varepsilon, \varepsilon)$ and thus
	$$
	\frac{f_{a_{\varepsilon}^*,b_{\varepsilon}}(u_{\varepsilon})}{\mathbf{d}(u_{\varepsilon}, \mathcal{S}(a_{\varepsilon}^*,b_{\varepsilon}))}=\frac{\varepsilon^2}{\sqrt{2}\varepsilon}
	=\frac{\varepsilon}{\sqrt{2}},
	$$
	which implies that  $\sigma(a_{\varepsilon}^*,b_{\varepsilon})<\varepsilon$.
\end{example}




\bibliography{WTY}

\bibliographystyle{unsrt}

\end{document}